\documentclass[10pt]{amsart}

\usepackage{amsmath}
\usepackage{amssymb}
\usepackage{euscript}
\usepackage[all]{xy}

\newtheorem{thm}{Theorem}
\newtheorem{prop}[thm]{Proposition}

\newcounter{exc}

\newcommand{\acat}[1]{\mathsf{a}_{#1}}

\newcommand{\sets}{\mathsf{Sets}}

\newcommand{\morcat}[1]{\mathsf{Mor \; {#1}}}
\newcommand{\mcat}[1]{\mathsf{Mod}(#1)}

\newcommand{\presh}{\mathsf{PreSh}}
\newcommand{\shc}[1]{\mathsf{Sh}(X,#1)}

\newcommand{\qcoh}[1]{\mathsf{QCoh}(#1)}
\newcommand{\compl}{\mathsf{Compl}}
\newcommand{\oc}[1]{\mathsf{#1}}

\newcommand{\sh}[1]{\mathcal{#1}}
\newcommand{\sa}{\mathcal{D}}
\newcommand{\sfam}{\mathcal{F}}

\newcommand{\defs}{\mathsf{Def}_{\sfam}}
\newcommand{\defm}[1]{\mathsf{Def}_{#1}}

\newcommand{\ch}{\mathsf{H}}
\newcommand{\cpd}{*}
\newcommand{\hc}{\mathsf{HC}}
\newcommand{\hh}{\mathsf{HH}}

\newcommand{\dc}{\mathsf{D}}

\newcommand{\reg}[1]{\mathcal{O}_{#1}}

\newcommand{\g}{\mathbf{g}}
\newcommand{\ue}{\textup{U}}

\newcommand{\pp}{\mathbf{P}}
\newcommand{\rst}[1]{|_{#1}}

\DeclareMathOperator{\coker}{coker}
\DeclareMathOperator{\der}{Der}
\DeclareMathOperator{\enm}{End}
\DeclareMathOperator{\ext}{Ext}

\DeclareMathOperator{\hmm}{Hom}
\DeclareMathOperator{\id}{id}

\title[Computing noncommutative deformations]{Computing noncommutative
global \\ deformations of D-modules}
\author{Eivind Eriksen}
\address{Oslo University College \\ Postboks 4 St. Olavs Plass \\
0130 Oslo, Norway}
\email{eeriksen@hio.no}

\begin{document}

\maketitle

\begin{abstract}
Let $(X,\sa)$ be a D-scheme in the sense of Beilinson and Bernstein,
given by an algebraic variety $X$ and a morphism $\reg X \to \sa$ of
sheaves of rings on $X$. We consider noncommutative deformations of
quasi-coherent sheaves of left $\sa$-modules on $X$, and show how to
compute their pro-representing hulls. As an application, we compute
the noncommutative deformations of the left $\sh D_X$-module $\reg
X$ when $X$ is any elliptic curve.
\end{abstract}

\section*{Introduction}

Let $k$ be an algebraically closed field of characteristic $0$, and let
$X$ be an algebraic variety over $k$, i.e. an integral separated scheme
of finite type over $k$. A \emph{D-algebra} in the sense of Beilinson,
Bernstein \cite{be-be93} is a sheaf $\sa$ of associative rings on $X$,
together with a morphism $i: \reg X \to \sa$ of sheaves of rings on
$X$, such that the following conditions hold: (i) $\sa$ is
quasi-coherent as a left and right $\reg X$-module via $i$, and (ii)
for any open subset $U \subseteq X$ and any section $P \in \sa(U)$,
there is an integer $n \ge 0$ such that
    \[ [ \dots [ [ P,a_1], a_2 ], \dots ,a_n] = 0 \]
for all sections $a_1, \dots, a_n \in \reg X(U)$, where $[P,Q] = PQ-QP$
is the commutator in $\sa(U)$. When $\sa$ is a D-algebra on $X$, the
ringed space $(X, \sa)$ is called a \emph{D-scheme}.

Let us denote the sheaf of $k$-linear differential operators on $X$ by
$\sh D_X$, and for any Lie algebroid $\g$ of $X/k$, let us denote the
enveloping D-algebra of $\g$ by $\ue(\g)$. We see that $\sh D_X$ and
$\ue(\g)$ are examples of noncommutative D-algebras on $X$, and that
$\reg X$ is an example of a commutative D-algebra on $X$.

Let us define a \emph{D-module} on a D-scheme $(X, \sa)$ to be a
quasi-coherent sheaf of left $\sa$-modules on $X$. In Eriksen
\cite{er06}, we developed a noncommutative global deformation theory
of D-modules that generalizes the usual (commutative) global
deformation theory of D-modules, and the noncommutative deformation
theory of modules in the affine case, due to Laudal. In section
\ref{s:defm} - \ref{s:comp-defm}, we review the essential parts of
this theory, including the global Hochschild cohomology and the
global obstruction calculus, all the time with a view towards
concrete computations.

The purpose of this paper is to show how to apply the theory in
order to compute noncommutative global deformations of interesting
D-modules. In section \ref{s:example}, we consider the
noncommutative deformation functor $\defm{\reg X}: \acat 1 \to
\sets$ of $\reg X$ considered as a left $\sh D_X$-module when $X$ is
any elliptic curve over $k$. Recall that in this case, a
quasi-coherent $\sh D_X$-module structure on $\reg X$ is the same as
an integrable connection on $\reg X$, and according to a theorem due
to Andr\'e Weil, see Weil \cite{weil38} and also Atiyah
\cite{ati57}, a line bundle admits an integrable connection if and
only if it has degree zero.

We show that the noncommutative deformation functor $\defm{\reg X}:
\acat 1 \to \sets$ has pro-representing hull $H =
k{\ll}t_1,t_2{\gg}/(t_1 t_2 - t_2 t_1) \cong k[[t_1, t_2]]$ that is
commutative, smooth and of dimension two. We also compute the
corresponding versal family in concrete terms, and remark that it
does not admit an algebraization.

\section{Noncommutative global deformations of D-modules} \label{s:defm}

Let $(X,\sa)$ be a D-scheme, and let $\qcoh \sa$ be the category of
quasi-coherent sheaves of left $\sa$-modules on $X$. This is the full
subcategory of $\shc \sa$, the category of sheaves of left
$\sa$-modules on $X$, consisting of quasi-coherent sheaves. We recall
that a sheaf $\sh F$ of left $\sa$-modules on $X$ is quasi-coherent if
for every point $x \in X$, there exists an open neighbourhood $U
\subseteq X$ of $x$, free sheaves $\sh L_0, \sh L_1$ of left $\sa \rst
U$-modules on $U$, and an exact sequence $0 \gets \sh F \rst U \gets
\sh L_0 \gets \sh L_1$ of sheaves of left $\sa \rst U$-modules on $U$.
We shall refer to the quasi-coherent sheaves of left $\sa$-modules on
$X$ as \emph{D-modules} on the D-scheme $(X, \sa)$.

For any D-scheme $(X,\sa)$, $\qcoh \sa$ is an Abelian $k$-category, and
we consider noncommutative deformations in $\qcoh \sa$. For any finite
family $\sfam = \{ \sh F_1, \dots, \sh F_p \}$ of quasi-coherent left
$\sa$-modules on $X$, there is a noncommutative deformation functor
$\defs^{qc}: \acat p \to \sets$ of $\sfam$ in $\qcoh \sa$, generalizing
the noncommutative deformation functor of modules introduced in Laudal
\cite{lau02}. We shall provide a brief description of $\defs^{qh}$
below; see Eriksen \cite{er06} for further details.

We recall that the objects of the category $\acat p$ of $p$-pointed
\emph{noncommutative Artin rings} are Artinian rings $R$, together
with structural morphisms $f: k^p \to R$ and $g: R \to k^p$, such
that $g \circ f = \id$ and the radical $I(R) = \ker(g)$ is
nilpotent. The morphisms are the natural commutative diagrams. For
any $R \in \acat p$, there are $p$ isomorphism classes of simple
left $R$-modules, represented by $\{ k_1, k_2, \dots, k_p \}$, where
$k_i = 0 \times \dots \times k \times \dots \times 0$ is the $i$'th
projection of $k^p$ for $1 \le i \le p$.

We remark that any $R \in \acat p$ is a $p \times p$ \emph{matrix
ring}, in the sense that there are $p$ indecomposable idempotents
$\{ e_1, \dots, e_p \}$ in $R$ with $e_1 + \dots + e_p = 1$ and a
decomposition $R = \oplus R_{ij}$, given by $R_{ij} = e_i R e_j$,
such that elements of $R$ multiply as matrices. We shall therefore
use matrix notation, and write $R = ( R_{ij} )$ when $R \in \acat
p$, and $( V_{ij} ) = \oplus V_{ij}$ when $\{ V_{ij}: 1 \le i,j \le
p \}$ is any family of vector spaces.

For any $R \in \acat p$, a lifting of $\sfam$ to $R$ is a sheaf $\sh
F_R$ in $\qcoh \sa$ with a compatible right $R$-module structure,
together with isomorphisms $\eta_i: \sh F_R \otimes_R k_i \to \sh F_i$
in $\qcoh \sa$ for $1 \le i \le p$, such that $\sh F_R(U) \cong ( \sh
F_i(U) \otimes_k R_{ij} )$ as right $R$-modules for all open subsets $U
\subseteq X$. We say that two liftings $(\sh F_R, \eta_i)$ and $(\sh
F'_R, \eta'_i)$ are equivalent if there is an isomorphism $\tau: \sh
F_R \to \sh F'_R$ of $\sa$-$R$ bimodules on $X$ such that $\eta'_i
\circ (\tau \otimes_k k_i) = \eta_i$ for $1 \le i \le p$, and denote
the set of equivalence classes of liftings of $\sfam$ to $R$ by
$\defs^{qc}(R)$. This defines the noncommutative deformation functor
$\defs^{qc}: \acat p \to \sets$.

\section{Computing noncommutative global deformations}
\label{s:comp-defm}

Let $(X, \sa)$ be a D-scheme, and let $\oc U$ be an open affine cover
of $X$ that is finite and closed under intersections. We shall explain
how to compute noncommutative deformations in $\qcoh \sa$ effectively
using the open cover $\oc U$.

We may consider $\oc U$ as a small category, where the objects are the
open subsets $U \in \oc U$, and the morphisms from $U$ to $V$ are the
(opposite) inclusions $U \supseteq V$. There is a natural forgetful
functor $\qcoh \sa \to \presh(\oc U, \sa)$, where $\presh(\oc U, \sa)$
is the Abelian $k$-category of (covariant) presheaves of left
$\sa$-modules on $\oc U$. For any finite family $\sfam$ in $\qcoh \sa$,
this forgetful functor induces an isomorphism of noncommutative
deformation functors $\defs^{qc} \to \defs$, where $\defs^{qc}: \acat p
\to \sets$ is the noncommutative deformation functor of $\sfam$ in
$\qcoh \sa$ defined in section \ref{s:defm}, and $\defs: \acat p \to
\sets$ is the noncommutative deformation functor of $\sfam$ in
$\presh(\oc U, \sa)$, defined in a similar way; see Eriksen \cite{er06}
for details.

\begin{thm} \label{t:defm-hull}
Let $(X,\sa)$ be a D-scheme, and let $\sfam = \{ \sh F_1, \dots, \sh
F_p \}$ be a finite family in $\qcoh \sa$. If the global Hochschild
cohomology $(\hh^n(\oc U, \sa, \hmm_k(\sh F_j, \sh F_i)))$ has finite
$k$-dimension for $n=1,2$, then the noncommutative deformation functor
$\defs^{qc}: \acat p \to \sets$ has a pro-representing hull $H =
H(\sfam)$, completely determined by $( \hh^n(\oc U, \sa, \hmm_k(\sh
F_j, \sh F_i)) )$ for $n=1,2$ and their generalized Massey products.
\end{thm}

In fact, there is a constructive proof of the fact that $\defs: \acat p
\to \sets$ of $\sfam$ in $\presh(\oc U, \sa)$ has a pro-representing
hull; see Eriksen \cite{er06} for details. The construction in this
proof uses the global Hochschild cohomology $( \hh^n(\oc U, \sa,
\hmm_k(\sh F_j, \sh F_i)) )$ for $n=1,2$, and the obstruction calculus
of $\defs$, which can be expressed in terms of generalized Massey
products on these cohomology groups. We give a brief description of the
global Hochschild cohomology and the obstruction calculus below.

\subsection{Cohomology}

For any presheaves $\sh F, \sh G$ of left $\sa$-modules on $\oc U$, we
recall the definition of the \emph{global Hochschild cohomology}
$\hh^n(\oc U, \sa, \hmm_k(\sh F, \sh G))$ of $\sa$ with values in the
bimodule $\hmm_k(\sh F, \sh G)$ on $\oc U$. For any (opposite)
inclusion $U \supseteq V$ in $\oc U$, we consider the Hochschild
complex $\hc^{\cpd}(\sa(U), \hmm_k(\sh F(U), \sh G(V))$ of $\sa(U)$
with values in the bimodule $\hmm_k(\sh F(U), \sh G(V))$. We define the
category $\morcat{\oc U}$ to have (opposite) inclusions $U \supseteq V$
in $\oc U$ as its objects, and nested inclusions $U' \supseteq U
\supseteq V \supseteq V'$ in $\oc U$ as its morphisms from $U \supseteq
V$ to $U' \supseteq V'$. It follows that we may consider the Hochschild
complex
    \[ \hc^{\cpd}(\sa, \hmm_k(\sh F, \sh G)): \morcat{\oc U} \to
    \compl(k) \]
as a functor on $\morcat{\oc U}$. The global Hochschild complex
$\hc^{\cpd}(\oc U, \sa, \hmm_k(\sh F, \sh G))$ is the total complex
of the double complex $\dc^{\cpd \cpd} = \dc^{\cpd}(\oc U,
\hc^{\cpd}(\sa, \hmm_k(\sh F, \sh G)))$, where $\dc^{\cpd}(\oc U,-):
\presh(\morcat{\oc U},k) \to \compl(k)$ is the resolving complex of
the projective limit functor; see Laudal \cite{lau79} for details.
Finally, we define the global Hochschild cohomology $\hh^n(\oc U,
\sa, \hmm_k(\sh F, \sh G))$ to be the cohomology of the global
Hochschild complex $\hc^{\cpd}(\oc U, \sa, \hmm_k(\sh F, \sh G))$.

We note that $\ch^n(\hc^{\cpd}(\sa(U), \hmm_k(\sh F(U), \sh G(V))))
\cong \ext^n_{\sa(U)}(\sh F(U), \sh G(V))$ for any $U \supseteq V$ in
$\oc U$ since $k$ is a field. Hence there is a spectral sequence
converging to the global Hochschild cohomology $\hh^n(\oc U, \sa,
\hmm_k(\sh F, \sh G))$ with
\begin{equation}
E^{pq}_2 = \ch^p(\oc U, \ext^q_{\sa}(\sh F, \sh G)),
\end{equation}
where $\ch^p(\oc U,-) = \ch^p(\dc^{\cpd}(\oc U,-))$ and we consider
$\ext^q_{\sa}(\sh F, \sh G): \morcat{\oc U} \to \mcat k$ as a functor
on $\morcat{\oc U}$, given by $\{ U \supseteq V \} \mapsto
\ext^q_{\sa(U)}(\sh F(U), \sh G(V))$ for all $q \ge 0$.

\subsection{Obstruction calculus}

Let $R \in \acat p$ and let $I = I(R)$ be the radical of $R$. For
any lifting $\sh F_R \in \defs(R)$ of the family $\sfam$ in
$\presh(\oc U, \sa)$ to $R$, we have that $\sh F_R(U) \cong ( \sh
F_i(U) \otimes_k R_{ij} )$ as a right $R$-module for all $U \in \oc
U$. Moreover, the lifting $\sh F_R$ is completely determined by the
left multiplication of $\sa(U)$ on $\sh F_R(U)$ for all $U \in \oc
U$ and the restriction map $\sh F_R(U) \to \sh F_R(V)$ for all $U
\supseteq V$ in $\oc U$. Let us write $Q^R(U,V) = ( \hmm_k(\sh
F_j(U), \sh F_i(V) \otimes_k R_{ij}))$ and $Q^R(U) = Q^R(U,U)$ for
all $U \supseteq V$ in $\oc U$. Then $\sh F_R \in \defs(R)$ is
completely described by the following data:
\begin{enumerate}
    \item \label{i:20}
    For all $U \in \oc U$, a $k$-algebra homomorphism $L(U):
    \sa(U) \to Q^R(U)$ satisfying $L(U)(P)(f_j) = Pf_j \otimes e_j +
    ( \sh F_i(U) \otimes_k I_{ij} )$ for all $P \in \sa(U), \; f_j \in
    \sh F_j(U)$.
    \item \label{i:11}
    For all inclusions $U \supseteq V$ in $\oc U$, a restriction
    map $L(U,V) \in Q^R(U,V)$ satisfying $L(U,V)(f_j) = (f_j \rst V)
    \otimes e_j + ( \sh F_i(V) \otimes I_{ij} )$ for all $f_j \in \sh
    F_j(U)$ and $L(U,V) \circ L(U)(P) = L(V)(P \rst V) \circ L(U,V)$
    for all $P \in \sa(U)$.
    \item \label{i:02}
    For all inclusions $U \supseteq V \supseteq W$ in $\oc U$,
    we have $L(V,W) L(U,V) = L(U,W)$ and $L(U,U) = \id$.
\end{enumerate}
A \emph{small surjection} in $\acat p$ is a surjective morphism $u: R
\to S$ in $\acat p$ such that $IK = KI = 0$, where $K = \ker(u)$ and $I
= I(R)$ is the radical of $R$. To describe the obstruction calculus of
$\defs$, it is enough to consider the following problem: Given a small
surjection $u: R \to S$ and a deformation $\sh F_S \in \defs(S)$, what
are the possible liftings of $\sh F_S$ to $R$? The answer is given by
the following proposition; see Eriksen \cite{er06} for details:

\begin{prop} \label{p:obstr}
Let $u: R \to S$ be a small surjection in $\acat p$ with kernel $K$,
and let $\sh F_S \in \defs(S)$ be a deformation. Then there exists a
canonical obstruction
    \[ o(u, \sh F_S) \in ( \hh^2(\oc U, \sa, \hmm_k(\sh F_j, \sh F_i))
    \otimes_k K_{ij} ) \]
such that $o(u, \sh F_S) = 0$ if and only if there exists a
deformation $\sh F_R \in \defs(R)$ lifting $\sh F_S$ to $R$.
Moreover, if $o(u, \sh F_S) = 0$, then there is a transitive and
effective action of $( \hh^1(\oc U, \sa, \hmm_k(\sh F_j, \sh F_i))
\otimes_k K_{ij} )$ on the set of liftings of $\sh F_S$ to $R$.
\end{prop}

In fact, let the deformation $\sh F_S \in \defs(S)$ be given by
$L^S(U): \sa(U) \to Q^S(U)$ and $L^S(U,V) \in Q^S(U,V)$ for all $U
\supseteq V$ in $\oc U$, and let $\sigma: S \to R$ be a $k$-linear
section of $u: R \to S$ such that $\sigma(e_i) = e_i$ and
$\sigma(S_{ij}) \subseteq R_{ij}$ for $1 \le i,j \le p$. We consider
$L^R(U): \sa(U) \to Q^R(U)$ given by $L^R(U) = \sigma \circ L^S(U)$
and $L^R(U,V) = \sigma(L^S(U,V))$ for all $U \supseteq V$ in $\oc
U$. The obstruction $o(U,\sh F_S)$ for lifting $\sh F_S$ to $R$ is
given by
\begin{enumerate}
    \item
    $(P,Q) \mapsto L^R(U)(PQ) - L^R(U)(P) \circ L^R(U)(Q)$ for all $U
    \in \oc U, \; P,Q \in \sa(U)$
    \item
    $P \mapsto L^R(U,V) \circ L^R(U)(P) - L^R(V)(P \rst V) \circ
    L^R(U,V)$ for all $U \supseteq V$ in $\oc U$, $P \in \sa(U)$
    \item
    $L^R(V,W) \circ L^R(U,V) - L^R(U,W)$ for all $U \supseteq V \supseteq
    W$ in $\oc U$
\end{enumerate}
We see that these expressions are exactly the obstructions for
$L^R(U)$ and $L^R(U,V)$ to satisfy conditions (\ref{i:20}) -
(\ref{i:02}) in the characterization of $\defs(R)$ given above.

\section{Calculations for D-modules on elliptic curves}
\label{s:example}

Let $X \subseteq \pp^2$ be the irreducible projective plane curve
given by the homogeneous equation $f = 0$, where $f = y^2 z - x^3 -
a x z^2 - b z^3$ for fixed parameters $(a,b) \in k^2$. We assume
that $\Delta = 4a^3 + 27b^2 \neq 0$, so that $X$ is smooth and
therefore an elliptic curve over $k$. We shall compute the
noncommutative deformations of $\reg X$, considered as a
quasi-coherent left $\sh D_X$-module via the natural left action of
$\sh D_X$ on $\reg X$.

We choose an open affine cover $\oc U = \{ U_1, U_2, U_3 \}$ of $X$
closed under intersections, given by $U_1 = D_+(y)$, $U_2 = D_+(z)$
and $U_3 = U_1 \cap U_2$. We recall that the open subset $D_+(h)
\subseteq X$ is given by $D_+(h) = \{ p \in X: h(p) \neq 0 \}$ for
$h = y$ or $h = z$. It follows from the results in section
\ref{s:comp-defm} that the noncommutative deformation functor
$\defm{\reg X}: \acat 1 \to \sets$ has a pro-representing hull $H$,
completely determined by the global Hochschild cohomology groups and
some generalized Massey products on them. We shall therefore compute
$\hh^n(\oc U, \sa, \enm_k(\reg X))$ for $n = 1,2$.

It is known that $\sh D_X(U)$ is a simple Noetherian ring of global
dimension one and that $\reg X(U)$ is a simple left $\sh
D_X(U)$-module for any open affine subset $U \subseteq X$; see for
instance Smith, Stafford \cite{sm-st88}. The functor $\ext^q_{\sh
D_X}(\reg X, \reg X): \morcat{ \oc U} \to \mcat k$ is therefore
given by $\ext^q_{\sh D_X}(\reg X, \reg X) = 0$ for $q \ge 2$ and
$\enm_{\sh D_X}(\reg X) = k$. Since the spectral sequence for global
Hochschild cohomology given in section \ref{s:comp-defm}
degenerates,
\begin{align*}
    \hh^n(\oc U, \sh D_X, \enm_k(\reg X)) & \cong \ch^{n-1}(\oc U,
    \ext^1_{\sh D_X}(\reg X, \reg X)) \text{ for } n \ge 1 \\
    \hh^0(\oc U, \sh D_X, \enm_k(\reg X)) & \cong k
\end{align*}
We shall compute $\ext^1_{\sh D_X}(\reg X, \reg X)$ and use this to
find $\ch^{n-1}(\oc U, \ext^1_{\sh D_X}(\reg X, \reg X))$ for $n =
1,2$.

Let $A_i = \reg X(U_i)$ and $D_i = \sh D_X(U_i)$ for $i=1,2,3$. We see
that $A_1 \cong k[x,z]/(f_1)$ and $A_2 \cong k[x,y]/(f_2)$, where $f_1
= z - x^3 - axz^2 - bz^3$ and $f_2 = y^2 - x^3 - ax - b$. Moreover, we
have that $\der_k(A_i) = A_i \partial_i$ and $D_i = A_i
\langle\partial_i\rangle$ for $i=1,2$, where
\begin{align*}
    \partial_1 & = (1 - 2axz - 3bz^2) \; \partial / \partial x +
    (3x^2 + az^2) \; \partial / \partial z \\
    \partial_2 & = -2y \; \partial / \partial x - (3x^2 + a) \;
    \partial / \partial y
\end{align*}
On the intersection $U_3 = U_1 \cap U_2$, we choose an isomorphism $A_3
\cong k[x,y,y^{-1}]/(f_3)$ with $f_3 = f_2$, and see that $\der_k(A_3)
= A_3 \partial_3$ and $D_3 = A_3 \langle\partial_3\rangle$ for
$\partial_3 =
\partial_2$. The restriction maps of $\reg X$ and $\sh D_X$, considered
as presheaves on $\oc U$, are given by
    \[ x \mapsto xy^{-1}, \; z \mapsto y^{-1}, \; \partial_1 \mapsto
    \partial_3 \]
for the inclusion $U_1 \supseteq U_3$, and the natural localization map
for $U_2 \supseteq U_3$. Finally, we find a free resolution of $A_i$ as
a left $D_i$-module for $i=1,2,3$, given by
\begin{align*}
0 \gets A_i & \gets D_i \xleftarrow{\cdot \partial_i} D_i \gets 0
\end{align*}
and use this to compute $\ext^1_{D_i}(A_i, A_j) \cong
\coker(\partial_i \rst{U_j}: A_j \to A_j)$ for all $U_i \supseteq
U_j$ in $\oc U$. We see that $\ext^1_{D_i}(A_i, A_3) \cong
\coker(\partial_3: A_3 \to A_3)$ is independent of $i$, and find the
following $k$-linear bases for $\ext^1_{D_i}(A_i,A_j)$:

\vspace{5mm}
\begin{center}
\begin{tabular}{|l|l|l|}
    \hline
    & $a \neq 0:$ & $a = 0:$ \\
    \hline
    $U_1 \supseteq U_1$ & $1, z, z^2, z^3$ & $1, z, x, xz$ \\
    $U_2 \supseteq U_2$ & $1, y^2$ & $1, x$ \\
    $U_3 \supseteq U_3$ & $x^2 y^{-1}, 1, y^{-1}, y^{-2}, y^{-3}$
    & $x^2 y^{-1}, 1, y^{-1}, x, x y^{-1}$ \\
    \hline
\end{tabular}
\end{center}
\vspace{5mm}

\noindent The functor $\ext^1_{\sh D_X}(\reg X, \reg X): \morcat{\oc U}
\to \mcat k$ defines the following diagram in $\mcat k$, where the maps
are induced by the restriction maps on $\reg X$:
\[
\xymatrix{ \ext^1_{D_1}(A_1,A_1) \ar[d] & & \ext^1_{D_2}(A_2,A_2)
\ar[d] \\
\ext^1_{D_1}(A_1,A_3) & \ext^1_{D_3}(A_3,A_3) \ar@{=}[l] \ar@{=}[r] &
\ext^1_{D_2}(A_2,A_3) }
\]
We use that $15 \; y^2 = \Delta \; y^{-2}$ in
$\ext^1_{D_3}(A_3,A_3)$ when $a \neq 0$ and that $-3b \; xy^{-2} =
x$ in $\ext^1_{D_3}(A_3,A_3)$ when $a = 0$ to describe these maps in
the given bases. We compute $\ch^{n-1}(\oc U, \ext^1_{\sh D_X}(\reg
X, \reg X))$ for $n = 1,2$ using the resolving complex
$\dc^{\cpd}(\oc U, -)$; see Laudal \cite{lau79} for definitions. In
particular, we see that $\ch^0(\oc U, \ext^1_{\sh D_X}(\reg X, \reg
X))$ consists of all pairs $(h_1, h_2)$ with $h_i \in
\ext^1_{D_i}(A_i,A_i)$ for $i = 1,2$ that satisfies the condition
$h_1 \rst{U_3} = h_2 \rst{U_3}$. Moreover, $\ch^1(\oc U, \ext^1_{\sh
D_X}(\reg X, \reg X))$ consists of all pairs $(h_{13}, h_{23})$ with
$h_{ij} \in \ext^1_{D_i}(A_i,A_j)$ for $(i,j) = (1,3), (2,3)$,
modulo the pairs of the form $(h_1 \rst{U_3} - h_3, h_2 \rst{U_3} -
h_3)$ for triples $(h_1, h_2, h_3)$ with $h_i \in
\ext^1_{D_i}(A_i,A_i)$ for $i = 1,2,3$. We find the following
$k$-linear bases:

\vspace{5mm}
\begin{center}
\begin{tabular}{|l|l|l|}
    \hline
    & $a \neq 0:$ & $a = 0:$ \\
    \hline
    $n=1$ & $\xi_1 = (1,1,1), \; \xi_2 = (\Delta z^2, 15 y^2, \Delta
    y^{-2})$ & $\xi_1 = (1,1,1), \; \xi_2 = (-3b \; xz, x, x)$ \\
    $n=2$ & $\omega = (0,0,0,0,6a x^2 y^{-1})$ & $\omega =
    (0,0,0,0,x^2 y^{-1})$ \\
    \hline
\end{tabular}
\end{center}
\vspace{5mm}

\noindent We recall that $\xi_1, \xi_2$ and $\omega$ are represented by
cocycles of degree $p = 0$ and $p = 1$ in the resolving complex
$\dc^{\cpd}(\oc U, \ext^1_{\sh D_X}(\reg X, \reg X))$, where
    \[ \dc^p(\oc U, \ext^1_{\sh D_X}(\reg X, \reg X)) = \prod_{U_0
    \supseteq \dots \supseteq U_p} \ext^1_{\sh D_X}(\reg X,\reg X)(U_0
    \supseteq U_p) \]
and the product is indexed by $\{ U_1 \supseteq U_1, U_2 \supseteq U_2,
U_3 \supseteq U_3 \}$ when $p = 0$, and $\{ U_1 \supseteq U_1, U_2
\supseteq U_2, U_3 \supseteq U_3, U_1 \supseteq U_3, U_2 \supseteq U_3
\}$ when $p = 1$.

This proves that the noncommutative deformation functor $\defm{\reg X}:
\acat 1 \to \sets$ of the left $\sh D_X$-module $\reg X$ has tangent
space $\hh^1(\oc U, \sh D_X, \enm_k(\reg X)) \cong k^2$ and obstruction
space $\hh^2(\oc U, \sh D_X, \enm_k(\reg X)) \cong k$ for any elliptic
curve $X$ over $k$, and a pro-representing hull $H =
k{\ll}t_1,t_2{\gg}/(F)$ for some noncommutative power series $F \in
k{\ll}t_1,t_2{\gg}$.

We shall compute the noncommutative power series $F$ and the versal
family $\sh F_H \in \defm{\reg X}(H)$ using the obstruction
calculus. We choose base vectors $t_1^*, t_2^*$ in $\hh^1(\oc U,
\sa, \hmm_k(\reg X, \reg X))$, and representatives $(\psi_l, \tau_l)
\in \dc^{01} \oplus \dc^{10}$ of $t_l^*$ for $l = 1,2$, where
$\dc^{pq} = \dc^p(\oc U, \hc^q(\sa, \enm_k(\reg X)))$. We may choose
$\psi_l(U_i)$ to be the derivation defined by
\begin{equation*}
\psi_l(U_i)(P_i) = \begin{cases} 0 & \text{if $P_i \in A_i$} \\
\xi_l(U_i) \cdot \id_{A_i} & \text{if $P_i = \partial_i$}
\end{cases}
\end{equation*}
for $l = 1,2$ and $i = 1,2,3$, and $\tau_l(U_i \supseteq U_j)$ to be
the multiplication operator in $\hmm_{A_i}(A_i,A_j) \cong A_j$ given by
$\tau_1 = 0$, $\tau_2(U_i \supseteq U_i) = 0$ for $i = 1,2,3$ and

\vspace{5mm}
\begin{center}
\begin{tabular}{|l|l|l|}
    \hline
    $a \neq 0:$ & $a = 0:$ \\
    \hline
    $\tau_2(U_1 \supseteq U_3) = 0$ & $\tau_2(U_1 \supseteq U_3) =
    x^2y^{-1}$ \\
    $\tau_2(U_2 \supseteq U_3) = -4a^2y^{-1} - 3xy + 9bxy^{-1} -
    6ax^2y^{-1}$ & $\tau_2(U_2 \supseteq U_3) = 0$ \\
    \hline
\end{tabular}
\end{center}
\vspace{5mm}

Let $\acat 1(n)$ be the full subcategory of $\acat 1$ consisting of
all $R$ such that $I(R)^n = 0$ for $n \ge 2$. The restriction of
$\defm{\reg X}: \acat 1 \to \sets$ to $\acat 1(2)$ is represented by
$(H_2, \sh F_{H_2})$, where $H_2 = k \langle t_1,t_2 \rangle
/(t_1,t_2)^2$ and the deformation $\sh F_{H_2} \in
\defm{\reg X}(H_2)$ is defined by $\sh F_{H_2}(U_i) = A_i \otimes_k
H_2$ as a right $H_2$-module for $i = 1,2,3$, with left $D_i$-module
structure given by
\begin{equation*}
P_i (m_i \otimes 1) = P_i(m_i) \otimes 1 + \psi_1(U_i)(P_i)(m_i)
\otimes t_1 + \psi_2(U_i)(P_i)(m_i) \otimes t_2
\end{equation*}
for $i = 1,2,3$ and for all $P_i \in D_i, \; m_i \in A_i$, and with
restriction map for the inclusion $U_i \supseteq U_j$ given by
\begin{equation*}
m_i \otimes 1 \mapsto m_i \rst{U_j} \otimes 1 + \tau_2(U_i \supseteq
U_j) \; m_i \rst{U_j} \otimes t_2
\end{equation*}
for $i=1,2, \; j = 3$ and for all $m_i \in A_i$.

Let us attempt to lift the family $\sh F_{H_2} \in \defm{\reg X}(H_2)$
to $R = k{\ll}t_1,t_2{\gg}/(t_1,t_2)^3$. We let $\sh F_R(U_i) = A_i
\otimes_k R$ as a right $R$-module for $i = 1,2,3$, with left
$D_i$-module structure given by
\begin{equation*}
P_i (m_i \otimes 1) = P_i(m_i) \otimes 1 + \psi_1(U_i)(P_i)(m_i)
\otimes t_1 + \psi_2(U_i)(P_i)(m_i) \otimes t_2
\end{equation*}
for $i = 1,2,3$ and for all $P_i \in D_i, \; m_i \in A_i$, and with
restriction map for the inclusion $U_i \supseteq U_j$ given by
\begin{equation*}
m_i \otimes 1 \mapsto m_i \rst{U_j} \otimes 1 + \tau_2(U_i \supseteq
U_j) \; m_i \rst{U_j} \otimes t_2 + \frac{\tau_2(U_i \supseteq
U_j)^2}{2} \; m_i \rst{U_j} \otimes t_2^2
\end{equation*}
for $i=1,2, \; j = 3$ and for all $m_i \in A_i$. We see that $\sh
F_R(U_i)$ is a left $\sh D_X(U_i)$-module for $i = 1,2,3$, and that
$t_1 t_2 - t_2 t_1 = 0$ is a necessary and sufficient condition for
$\sh D_X$-linearity of the restriction maps for the inclusions $U_1
\supseteq U_3$ and $U_2 \supseteq U_3$. This implies that $\sh F_R$
is not a lifting of $\sh F_{H_2}$ to $R$. But if we define the
quotient $H_3 = R/(t_1 t_2 - t_2 t_1)$, we see that the family $\sh
F_{H_3} \in \defm{\reg X}(H_3)$ induced by $\sh F_R$ is a lifting of
$\sh F_{H_2}$ to $H_3$.

In fact, we claim that the restriction of $\defm{\reg X}: \acat 1
\to \sets$ to $\acat 1(3)$ is represented by $(H_3, \sh F_{H_3})$.
One way to prove this is to show that it is not possible to find any
lifting $\sh F'_R \in \defm{\reg X}(R)$ of $\sh F_{H_2}$ to $R$.
Another approach is to calculate the cup products $< t^*_i, t^*_j>$
in global Hochschild cohomology for $i,j = 1,2$, and this gives
\begin{align*}
    < t^*_1, t^*_2> = o^*, \; & < t^*_2, t^*_1> = - o^* &
    \text{ for } a \neq 0 \\
    < t^*_1, t^*_2> = o^*, \; & < t^*_2, t^*_1> = - o^* &
    \text{ for } a = 0
\end{align*}
where $o^* \in \hh^2(\oc U, \sh D_X, \enm_k(\reg X))$ is the base
vector corresponding to $\omega$. Since all other cup products vanish,
this implies that $F = t_1 t_2 - t_2 t_1 + (t_1,t_2)^3$.

Let $H = k{\ll}t_1,t_2{\gg}/(t_1t_2 - t_2 t_1)$. We shall show that
it is possible to find a lifting $\sh F_H \in \defm{\reg X}(H)$ of
$\sh F_{H_3}$ to $H$. We let $\sh F_H(U_i) = A_i \widehat{\otimes}_k
H$ as a right $H$-module for $i = 1,2,3$, with left $D_i$-module
structure given by
\begin{equation*}
P_i (m_i \otimes 1) = P_i(m_i) \otimes 1 + \psi_1(U_i)(P_i)(m_i)
\otimes t_1 + \psi_2(U_i)(P_i)(m_i) \otimes t_2
\end{equation*}
for $i = 1,2,3$ and for all $P_i \in \sh D_i, \; m_i \in A_i$, and
with restriction map for the inclusion $U_i \supseteq U_j$ given by
\begin{equation*}
    m_i \otimes 1 \mapsto \sum_{n=0}^{\infty} \frac{\tau_2(U_i
    \supseteq U_j)^n}{n!} \; m_1 \rst{U_j} \otimes t_2^n
    = \exp(\tau_2(U_i \supseteq U_j) \otimes t_2) \cdot (m_1
    \rst{U_j} \otimes 1)
\end{equation*}
for $i=1,2, \; j = 3$ and for all $m_i \in A_i$. This implies that
$(H, \sh F_H)$ is the pro-representing hull of $\defm{\reg X}$, and
that $F = t_1 t_2 - t_2 t_1$. We remark that the versal family $\sh
F_H$ does not admit an algebraization, i.e. an algebra
$H_{\text{alg}}$ of finite type over $k$ such that $H$ is a
completion of $H_{\text{alg}}$, together with a deformation in
$\defm{\reg X}(H_{\text{alg}})$ that induces the versal family $\sh
F_H \in \defm{\reg X}(H)$.

Finally, we mention that there is an algorithm for computing the
pro-representing hull $H$ and the versal family $\sh F_H$ using the
cup products and higher generalized Massey products on global
Hochschild cohomology. We shall describe this algorithm in a
forthcoming paper. In many situations, it is necessary to use the
full power of this machinery to compute noncommutative deformation
functors.

\bibliographystyle{amsplain}
\bibliography{main}

\providecommand{\bysame}{\leavevmode\hbox to3em{\hrulefill}\thinspace}
\providecommand{\MR}{\relax\ifhmode\unskip\space\fi MR }
\providecommand{\MRhref}[2]{%
  \href{http://www.ams.org/mathscinet-getitem?mr=#1}{#2}
}
\providecommand{\href}[2]{#2}
\begin{thebibliography}{1}

\bibitem{ati57}
M.~F. Atiyah, \emph{Complex analytic connections in fibre bundles}, Trans.
  Amer. Math. Soc. (1957), 181--207.

\bibitem{be-be93}
A.~Beilinson and J.~Bernstein, \emph{A proof of the {J}antzen conjectures},
  Adv. Soviet Math. \textbf{16} (1993), no.~1, 1--50.

\bibitem{er06}
E.~Eriksen, \emph{Noncommutative deformations of sheaves and presheaves of
  modules}, ArXiv: math.AG/0405234, 2005.

\bibitem{lau79}
O.~A. Laudal, \emph{Formal moduli of algebraic structures}, Lecture Notes in
  Mathematics, no. 754, Springer-Verlag, 1979.

\bibitem{lau02}
\bysame, \emph{Noncommutative deformations of modules}, Homology Homotopy Appl.
  \textbf{4} (2002), no.~2, part 2, 357--396.

\bibitem{sm-st88}
S.~Paul Smith and J.~T. Stafford, \emph{Differential operators on an affine
  curve}, Proc. London Math. Soc. (3) \textbf{56} (1988), no.~2, 229--259.

\bibitem{weil38}
A.~Weil, \emph{G\'en\'eralisation des fonctions ab\'eliennes}, J. Math. Pures
  Appl. \textbf{17} (1938), 47--87.

\end{thebibliography}

\end{document}